\documentclass[12pt]{article}
\def\dateref{October 2020}
\usepackage{amssymb, amsmath, amsthm, amscd, graphicx, dsfont, mathabx}
\usepackage[T1]{fontenc}
\newtheorem{lem}{Lemma}[section]

\newtheorem{thm}{Theorem}[section]

\newtheorem{defn}{Definition}[section]

\def\cB{{\mathcal B}}
\def\cE{{\mathcal E}}
\def\cG{{\mathcal G}}
\def\cL{{\mathcal L}}
\def\cP{{\mathcal P}}
\def\cR{{\mathcal R}}
\def\cV{{\mathcal V}}

\def\N{\mathbb{N}}
\def\Z{\mathbb{Z}}
\def\R{\mathbb{R}}

\def\tA{{\tt A}}
\def\tb{{\tt b}}

\def\colon{\,{:}\;}
\def\prf{\medbreak\noindent{\bf Proof.}\enspace}
\def\qed{\hspace*{\fill}\hbox{\vrule height 7pt \kern-.3pt
     \vbox{\hrule width 7pt
     \kern6.6pt\hrule width 7pt }\kern-.3pt\vrule height 7pt
     }\par}
\def\ra{\rightarrow}

\def\!{\mskip-\thinmuskip}

\begin{document}

\title{\bf Equilibrium measures of the natural extension of $\beta$-shifts}

\author{C.-E. Pfister\footnote{E-mail: charles.pfister@epfl.ch}\\
Section of Mathematics,
Faculty of Basic Sciences, EPFL\\
CH-1015 Lausanne, Switzerland
\and W.G. Sullivan\footnote{E-mail:
Wayne.Sullivan@ucd.ie}\\
               Department of Mathematics, UCD,\\
               Belfield, Dublin 4, Ireland}

 \date{\dateref}

\maketitle

\noindent
{\bf Abstract\colon}
We give a necessary and sufficient condition on $\beta$ of the natural extension of a  $\beta$-shift, so that any equilibrium measure 
for a function $\varphi$ of bounded total  oscillations
is a weak Gibbs measure.

\section{Introduction}\label{introduction}

We study equilibrium measures of the natural extension of  $\beta$-shifts. This is an interesting class of dynamical systems which have been studied in ergodic and number theory since the fundamental papers 
\cite{Re}, \cite{Pa}.
We want to determinate whether an equilibrium measure for a continuous function $\varphi$ is  a weak Gibbs measure. 
In \cite{PS3} we developed a method  based on a decoupling property, definition 2.3 in \cite{PS3}, which is a slightly weaker condition than condition $(D)$ in \cite{Ru} section \S4.1. The results
of \cite{PS3} are valid only for $\beta's$  such that the $\beta$-shift has the specification property and $\varphi$ is of bounded total oscillations. This  set  of 
$\beta's$  is the set $C_3$ in \cite{Sc}.
Schmeling proved that $C_3$  has Lebesgue measure $0$, but Hausdorff dimension $1$. 
For the more restricted class of  functions  verifying  Bowen condition \cite{Bo}, one has 
a stronger result. Under expansiveness and specification Haydn and Ruelle \cite{HR} proved the equivalence of equilibrium measures and Gibbs measures (in the sense of Bowen  \cite{Bo} and Capocaccia \cite{Ca}). 

However, using  basic ideas of  \cite{PS3} and \cite{CTY},
 it is possible to obtain a necessary and sufficient condition on $\beta$ such that for any function of  bounded total  oscillations all equilibrium measures are weak Gibbs measures.

To formulate precisely our main result, theorem \ref{thmmain}, we need to recall first some basic properties of 
$\beta$-shifts. This is done in subsection \ref{betashift}. In subsection \ref{subsectionboundedoscillations} we consider the class of functions of bounded total oscillations following section 3 of \cite{PS3} and
in subsection \ref{subsectionpressure} we consider the pressure, establishing two basic estimates for the
proof of theorem \ref{thmmain}.
Our main result is stated in subsection \ref{subsectionresult} and proved in section \ref{sectionproof}. We discuss briefly large deviations for empirical measures in section \ref{subsectionresult}.

\section{Setting and main result}

\subsection{Beta-shift}\label{betashift}

Let $\beta>1$ be fixed. The case $\beta\in\N$ is special and corresponds to the full shift.  From now on we assume that 
$\beta\not\in\N$.
For $t\in\R$, let
$ \lceil t\rceil :=\min\{i\in\Z\colon i\geq t\}$.
We define $\tb:=\lceil \beta\rceil$.  Consider the
$\beta$-expansion of $1$,
$$
1=\sum_{i=1}^\infty c_i\beta^{-i}\,,
$$
which is  given by the algorithm
$$
  r_0:=1,\; c_{i+1}:=\lceil \beta\, r_i\rceil-1,\;
  r_{i+1}:=\beta \,r_i-c_{i+1},\;i\in\Z_+\,,
$$
which insures that $r_i>0$ for all $i\in\Z_+$. It follows that $c_1=\lceil\beta\rceil-1>0$ and
$c^\beta:=(c_1,c_2,\ldots)$ cannot end with zeros only. For sequences
$(a_1,a_2,\ldots)$ and $(b_1,b_2,\ldots)$ the lexicographical
order is defined by $(a_1,a_2,\ldots)\prec (b_1,b_2,\ldots)$ if and only if for the
least index $i$ with $a_i\neq b_i$, $a_i<b_i$. Let
$\tA:=\{0,\ldots,\tb-1\}$; the (one-sided) $\beta$-shift is
$$
  X^\beta:=\left\{ x=(x_1,x_2,\ldots)\colon  x_i\in\tA,\;
    T^kx \preceq c^\beta\;\forall \,k\in\Z_+\right\}\,,
 $$
where $T$ is the left shift operator.
In particular $T^kc^\beta \preceq c^\beta$ for all $k\in\Z_+$,  so that $X^\beta$ is
a shift-invariant closed subset of $\tA^\N$ (with product topology). The language of the shift $X^\beta$ is denoted by $\cL^\beta$ and the set of the words of length $n$ by
$\cL^\beta_n$; the empty-word is $\epsilon$, $\cL^\beta_0=\{\epsilon\}$. A word is written $w_1\cdots w_n$ or simply $w$. The length of a word $w$ is written $|w|$.

The shift-space $X^\beta$  can be described by a labeled graph $\cG^\beta=(\cV,\cE^\beta)$ where $\cV:=\{q_j\colon j\in \Z_+\}$. The root of the graph is the vertex $q_0$.
There is an edge $q_0\ra q_0$, labeled by $k$, for each $k=0,\ldots,\tb-2$, and there is an edge $q_{j-1}\ra q_ j$ labeled by $c_j$  for each $j\in\N$.
Moreover, if the label $c_j$ of $q_{j-1}\ra q_j$ is different from $0$, then there are $c_j-1$ edges $q_{j-1}\ra q_0$ labeled by $0,\ldots, c_j-1$.
Each word $w_1\cdots w_n\in\cL^\beta$ can always be presented by a path of length $n$ in $\cG^\beta$ starting with vertex $q_0$. For a word $w\in\cL^\beta$ 
we define $q(w)$ as the end vertex of this path starting at $q_0$ and presenting $w$.
One can concatenate two words  $w$ and $w^\prime$ if and only if there is a path $\eta$ in $\cG^\beta$ presenting $w$ and a path $\eta^\prime$ in 
$\cG^\beta$ presenting $w^\prime$, so that $\eta$ ends at vertex $q$ and $\eta^\prime$ starts at vertex $q$. In particular
 one can concatenate $w$ with any words of $\cL^\beta$ if $q(w)=q_0$.
 There is a unique labeled path presenting the infinite sequence $c^\beta$, which is the path  $(q_0,q_1,q_2,\ldots)$. 
Let $\cP^\beta$ be the set  of the
prefixes of the sequence $c^\beta$, including the empty-word  $\epsilon$. 
Let $c_1\cdots c_n\in\cP^\beta$  and suppose that $c_{n+1}=c_{m+2}=\cdots c_{n+m}=0$, $c_{n+m+1}\not=0$. The word  $c_1\cdots c_n$ is presented by a path starting at $q_0$ and ending at $q_n$. Since there is only one out-going edge from each of the vertices $q_n,\ldots, q_{n+m-1}$, 
the only words $w$ with prefix $c_1\cdots c_n$ are the words 
$$
c_1\cdots c_n,\,
c_1\cdots c_n0,\, c_1\cdots c_n00,\,\ldots ,\,c_1\cdots c_n\underbrace{0\cdots 0}_{m},\, 
c_1\cdots c_n\underbrace{0\cdots 0}_{m}w^\prime\,, c_1\cdots c_k\,,
$$
where $w^\prime$ is any word of $\cL^\beta$ with first letter $0\leq w^\prime_1\leq c_{n+m+1}-1$, and $c_1\cdots c_k\in\cP^\beta$ with $k>n+m$.
For $u\in\cP^\beta$ we set
\begin{equation}\label{eq1}
z^\beta(u):=\begin{cases}
0 &\text{if $u=\epsilon$,}\\
p& \text{if $u=c_1\cdots c_\ell$, $c_{\ell+1}=\cdots =c_{\ell+p}=0$ and $c_{\ell+p+1}>0$.}
\end{cases}
\end{equation}
$z^\beta(u)$ is a measure of the obstruction to go from vertex $q(u)$ to vertex $q_0$.
We set
\begin{equation}\label{eq2}
\overline{z}^\beta(n):=\max\{z^{\beta}(u)\colon \text{$u$ prefix of $c^\beta$, $|u|\leq n$}\}\,.
\end{equation}
For each prefix $u=c_1\cdots c_n$ of $c^\beta$ we define a new word $\widehat{u}$ as follows.
Let $c_j$ be the last letter in $c_1\cdots c_n$ which is different from $0$. We set 
\begin{equation}\label{eq3}
\widehat{c}_\ell:=\begin{cases}
c_\ell & \text{if $\ell\not=j$}\\
c_j-1 &\text{if $\ell=j$.}
\end{cases}
\end{equation}
The word $\widehat{u}:=\widehat{c}_1\cdots\widehat{c}_n$ differs from $u=c_1\cdots c_n$ by a single letter and $q(\widehat{u})=q_0$. 
For any word $w\in\cL^\beta$ there is a unique decomposition of $w$ into 
\begin{equation}\label{eqs(v)}
w=vs(w)\quad\text{where $s(w)$ is the largest suffix of $w$ belonging to $\cP^\beta$.}
\end{equation}
We  extend  the definition \eqref{eq1} to any word $w$, by setting
$$
z^\beta(w):=z^\beta(s(w))\,,
$$
and extend
the transformation $u\mapsto\widehat{u}$ to any word by setting
\begin{equation}\label{eq4}
\widehat{w}:=\begin{cases}
w &\text{if $s(w)=\epsilon$}\\
v\widehat{u} &\text{if $s(w)=u$.}
\end{cases}
\end{equation}
By convention we set $\widehat{\epsilon}=\epsilon$.
The words $\widehat{w}$ can be freely concatenated since $q(\widehat{w})=q_0$ (see lemma \ref{lem1}).

\begin{lem}\label{lem0}
Let $a=a_1\cdots a_k$ and $b=b_1\cdots b_\ell$ be two prefixes of $c^\beta$. If $ab\in\cL^\beta$, then $ab$ is a prefix
of $c^\beta$.
\end{lem}

\prf
By hypothesis $a_i=c_i$, $i=1,\ldots,k$ and $b_j=c_j$, $j=1,\ldots, \ell$. Let $w=ab\in\cL^\beta$. Suppose that $w$ is not a prefix of $c^\beta$.  
Then there exists
$j$, $1\leq j\leq\ell$,  so that 
$$
c_m=w_{k+m}=c_{k+m}\quad  \text{ for $1\leq m<j$}\quad\text{and}\quad c_j=w_{k+j}<c_{k+j}\,.
$$
Hence,  $w_{k+1}\cdots w_{k+\ell}=c_1\cdots c_\ell\prec c_{k+1}\cdots c_{k+\ell}$, a contradiction with 
 $T^kc^\beta\preceq c^\beta$.
 \qed

\begin{lem}\label{lem1}
a) Let $w=vu$, $s(w)=u$. Then $q(v)=q_0$ and $s(\widehat{w})=s(v\widehat{u})=\epsilon$.\\
b) Let $p_1:=z^\beta(c_1)$, $c_1$ the first character of $c^\beta$. Then the mapping on  $\cL^\beta$, $w\mapsto\widehat{w}$, is at most
$(p_1+2)$-to-$1$, and $s(\widehat{w})=\epsilon$.
\end{lem}

\prf
a) Let $w=vu$, $u=s(w)$ and $\widehat{w}=v\widehat{u}$. If $q(v)\not=q_0$, then $s(v)\not=\epsilon$, so that by lemma \ref{lem0} $u$ is not maximal, a contradiction. Since $s(v)=\epsilon$ and $q(\widehat{u})=q_0$, $s(\widehat{w})=s(\widehat{u})=\epsilon$.

b) Let $\widehat{\cP}^\beta:=\{\widehat{u}\colon u\in\cP^\beta\}$. The list is ordered according to increasing length. The first $p_1+2$ 
words in the list $\widehat{\cP}^\beta$ are 
\begin{equation}\label{eqlist}
\epsilon,\,\widehat{c}_1, \,\widehat{c}_10, \,\widehat{c}_100,\ldots\,, \widehat{c}_1\underbrace{0\cdots0}_{p_1}\,.
\end{equation}
On the other hand, if $|\widehat{u}|>p_1+1$, then the first character of $\widehat{u}$ is $c_1$.

Let $\widehat{w}=w^\prime$. Let $\widehat{s}(w^\prime)$ be the largest suffix of $w^\prime$ among the first $|w^\prime|+1$ elements of the list 
$\widehat{\cP}^\beta$.  We write $w^\prime=v^\prime u^\prime$ with $u^\prime=\widehat{s}(w^\prime)$.
Let $w=v u$, $u=s(w)$, such that $\widehat{w}=w^\prime$. 
We have $|u|\leq |u^\prime|$, otherwise $\widehat{w}=w^\prime$ would imply that $\widehat{s}(w^\prime)$ is not maximal.
In particular if   $\widehat{s}(w^\prime)=\epsilon$, then $s(w)=\epsilon$ and $w^\prime=w$.  

Suppose that $\widehat{w}=w^\prime$ and $p_1+2\leq |u|<|u^\prime|$. Then  the first character of $u^\prime$ is $c_1$
and also the first character of $\widehat{u}$ is $c_1$ since $ |u|\geq p_1+2$. By hypothesis $|u|<|u^\prime|$. This implies that 
$v=v^\prime a$ with $a$ a prefix of $u^\prime$. Hence the first letter of $a$ is the first letter of $u^\prime$, which is $c_1$, and the letter
following $a$ is the first letter of $\widehat{u}$, which is $c_1$. We have $u^\prime=a\widehat{u}=\widehat{s}(w^\prime)\in \widehat{\cP}^\beta$.
By definition of the map $w\mapsto\widehat{w}$ (see
\eqref{eq4} and \eqref{eq3}) we conclude that $a$ is a prefix of $c^\beta$. By lemma \ref{lem0} $au$ is a prefix of $c^\beta$,
contradicting the maximality of $u$. Therefore $|u|=|u^\prime|$ and in this case the mapping  $w\mapsto\widehat{w}$ is $2$-to-$1$. 
In the remaining cases  $|s(w)|\leq p_1+1$.
Therefore the mapping $w\mapsto\widehat{w}$  is at most $(p_1+2)$-to-$1$.
\qed

\begin{defn}
The   the natural extension $\Sigma^\beta$ of $X^\beta$ is
$$
\Sigma^\beta=\{x\in\tA^\Z\colon \forall k\in\Z\,,\;(x_k,x_{k+1},\ldots)\in X^\beta\}\,.
$$
It is called hereafter simply \emph{$\beta$-shift}.
\end{defn}

The language of $\Sigma^\beta$ is also $\cL^\beta$.
Let $k<\ell$, $[k,\ell]=\{k,k+1,\ldots,\ell-1,\ell\}$, and $x\in \Sigma^\beta$. The projection $J_{[k,\ell]}\colon \Sigma^\beta\ra\cL^\beta$ is defined as
$$
x\mapsto J_{[k,\ell]}(x):=x_{[k,\ell]}\equiv x_kx_{k+1}\cdots x_{\ell}\,.
$$
Let $w=x_kx_{k+1}\cdots x_{k+m-1}\in\cL_m^\beta$.
We can always extend 
$w$ to the left by $0$, that is, there is $y\in\Sigma^\beta$, 
$y_j=0$, $j<k$ and $y_j=x_j$, $j=k,\ldots,k+m-1$.
 We can also extend $w$ to the right by $0$. If $q(w)=q_0$, this is clear. If $w=vu$,
$s(w)=u\not =\epsilon$, then
$u=c_1\cdots c_p$ for some $p\geq 1$. When $c_{p+1}\not=0$, we may change $c_{p+1}$ into $0$. When $c_{p+1}=\cdots=c_{p+r}=0$, but
$c_{p+r+1}\not=0$, we may change $c_{p+r+1}$ into $0$. Hence there exists $y\in\Sigma^\beta$, 
$y_j=x_j$, for all $j=k,\ldots,k+m-1$, and  
 $y_j=0$ for all $j<k$ and $j\geq m$.

\subsection{Functions of  bounded total  oscillations}\label{subsectionboundedoscillations}

We recall the definition of a function of  bounded total  oscillations. For details we refer to section 3 in \cite{PS3}. 
Let $\bar f\in C(\tA^\Z)$ be a continuous function defined on the \emph{full shift} $\tA^\Z$. On $\tA^\Z$ we define for each $i\in\Z$
$$
 \delta_i(\bar f):=\sup\{|\bar f(x)-\bar f(y)|\colon x,y\in \tA^\Z,\;
   \text{where $x_k=y_k$ for all $k\not=i$}\}
$$
and
$$
  \|\bar f\|_\delta:=\sum_{i\in \Z}\delta_i(\bar f)\,.
$$
A function has \emph{bounded total oscillations} if  $\|\bar f\|_\delta<\infty$.
On a subshift $X\subset \tA^\Z$  $\delta_i(f)$ may not make sense. If $f\in C(X)$ is a continuous function on $X$ and has a continuous extension
$\bar f\in C(\tA^\Z)$ on $\tA^Z$, we write
$$
f\approx \bar f\iff f(x)=\bar f(x)\quad\forall x\in X\,.
$$
Extension $\bar f$ of $f$ exists (proposition 3.2 \cite{PS3}). For $f\in C(X)$ we define
$$
\|f\|_\delta:=\inf\{\|\bar f\|_\delta\colon \bar f\in C(\tA^\Z)\,,\;\bar f\approx f\}\,.
$$
A function $f\in C(X)$ has \emph{bounded total oscillations} if $\|f\|_\delta<\infty$. Examples of functions of bounded total oscillations are given in \cite{PS3}. The set of bounded total oscillations is a Banach space $\cB(\Sigma^\beta)$ with the norm (proposition 3.1 \cite{PS3})
$$
 \| f\|_\delta + \sup_{x\in X} |f(x)|\,.
$$

We prove two basic estimates for functions with bounded total oscillations.
For convenience, from now on we write  $\bar f$ a continuous  extension  of $f$ to $\tA^\Z$.
The arguments do not require that $\bar f$  satisfies
$\|\bar f\|_\delta=\| f\|_\delta$ but just that $f\approx \bar f$ and  $\|\bar f\|_\delta<\infty$.
Fundamental to many of the arguments is the following lemma. 

\begin{lem}\label{lemfund}
Let $x,y\in X$ and $\Gamma:=\{j\colon x_j\not=y_j\}$. Then for $\Lambda\subset \Z$ and $f\in C(X)$,
$$
   \sum_{i\in\Lambda}|f(T^ix)- f(T^iy)|\leq
   \sum_{i\in\Lambda}\sum_{j\in\Gamma}\delta_{j-i} (\bar f)\leq\infty\,.
$$
\end{lem}
\prf
Since $\Gamma$ is at most countable, we can list the
elements of $\Gamma$, so that $\Gamma=\{j_1,j_2,\ldots\}$.  We define  a sequence of elements
of $\tA^\Z$ as follows. Let $z^{j_0}:=x$. For $j_\ell\in\Gamma$, set
$$
z^{j_\ell}_k:=
\begin{cases}
y_k &\text{if $k=j_1,\ldots,j_{\ell}$}\\
x_k &\text{if $k\in\Z\setminus\{j_1,\ldots,j_{\ell}\}$.}
\end{cases}
$$
Then
$$
|f(T^ix)- f(T^iy)|\leq \sum_{j_k\in\Gamma}|\bar f(T^iz^{j_{k-1}})-\bar f(T^iz^{j_k})|\leq 
\sum_{j_k\in\Gamma}\delta_{j_k}(\bar f\circ T^i)\,.
$$
The lemma follows from the identity
$$
\delta_j(\bar f\circ T^i)=\delta_{j-i}(\bar f)\,.
$$
Indeed
\begin{eqnarray*}
\delta_j(\bar f\circ T^i)
&=&
\sup\{|(\bar f\circ T^i)(x)-(\bar f\circ T^i)(y)|\colon x_k=y_k\;\forall k\not=j\}\\
&=&
\sup\{|\bar f(T^ix)-\bar f(T^iy)|\colon (T^{-i}T^ix)_k=(T^{-i}T^iy)_k\;\forall k\not=j\}\\
&=&
\sup\{|\bar f(x^\prime)-\bar f(y^\prime)|\colon (T^{-i}x^\prime)_k=(T^{-i}y^\prime)_k\;\forall k\not=j\}
=\delta_{j-i}(\bar f)\,.
\end{eqnarray*}
\qed

\begin{lem}\label{lem2}
Let $f$ be a function of bounded total oscillations on $X$.  Given $\varepsilon>0$, there exists 
$N_\varepsilon$ such that for $m\geq N_\varepsilon$
\begin{equation}\label{eq6} 
\sup\big\{\sum_{1\leq i\leq m}\big|f(T^ix)-f(T^iy)\big|\colon \text{$x,y$, 
$x_k=y_k$ $\forall k \in \{1,\ldots,m\}$}\big\}\leq m\varepsilon
\end{equation}
and
\begin{equation}\label{eq7} 
\sup\big\{\sum_{j\not\in \{1,\ldots,m\}}\big|f(T^jx)-f(T^jy)\big|\colon \text{$x,y$,  $x_k=y_k$  
$\forall k \not\in\{1,\ldots,m\}$}\big\}\leq m\varepsilon\,.
\end{equation}
\end{lem}

\prf
Let $\varepsilon>0$ be given.
There exists $r_{\varepsilon}$ so that
$\sum_{k:|k|>r_{\varepsilon}}\delta_k(\bar f)\leq \varepsilon/2$. If $m>2r_{\varepsilon}$,
then for $x_{[1,m]}=y_{[1,m]}$ the sum over $[1,m]$ of
$|T^i f(x)-T^i f(y)|$ can be written as over 
$$ 
  [1,r_\varepsilon]\cup[r_\varepsilon+1,m-r_\varepsilon]\cup[m+1-r_\varepsilon,m].
$$
For $i$ in the middle interval and $j\notin[1,m]$ we have $|i-j|>r_\varepsilon$, so that by lemma \ref{lemfund}
$$
| f(T^ix)- f(T^iy)|\leq \sum_{j\not\in\{1,\ldots,m\}}\delta_j(\bar f\circ T^i)\leq \sum_{k:|k|>r_{\varepsilon}}\delta_k(\bar f)\leq\varepsilon/2\,.
$$
For the $i$ in the outside intervals
we use $| f(T^ix)- f(T^iy)|\leq 2\| f\|$ to yield 
$$
\frac{1}{m}\sum_{1\leq i\leq m}\big|f(T^ix)-f(T^iy)\big|\leq 
\frac{4r_{\varepsilon}\|f\|}{m}+\frac{\varepsilon(m-2r_\varepsilon)}{2m}\leq\varepsilon m\quad\text{for $m$ large enough.}
$$
The proof of the second statement is similar.
\qed

\subsection{Equilibrium measure and pressure}\label{subsectionpressure}

In our setting   a shift-invariant (Borel) probability measure $\nu$ is an equilibrium measure for a continuous function $\varphi$ if and only if $\nu$ is a tangent functional to the pressure $p$ at $\varphi$ (see \cite{Wa}  theorems 8.2 and 9.5).

\begin{defn}\label{defntangentfunctional}
An invariant   probability measure $\nu$ is a \emph{tangent functional to the pressure $p$ at $\varphi$} if 
$$
p(\varphi+f)\geq p(\varphi)+\int f\,d\nu\,,\;\text{for all continuous functions $f$}\,.
$$
The set of tangent functionals to the pressure at $\varphi$ is denoted $\partial p(\varphi)$.
\end{defn}

For each  $n\in\N$ we choose a set $E^n$ with the following properties: 
\begin{equation}\label{eq8}
J_{[-n,n]}(E^n)=\cL^\beta_{2n+1}\quad\text{and}\quad
\big(x,x^\prime\in E^n\,,\;J_{[-n,n]}(x)=J_{[-n,n]}(x^\prime)\big)\implies x=x^\prime\,.
\end{equation}
Let  $\varphi$ be a continuous function and set
\begin{equation}\label{eq9}
\Xi^n(\varphi):=\sum_{x\in E^n}\exp\big(\sum_{j\in[-n,n]}\varphi(T^jx)\big)\quad\text{and}\quad P_{E^n}(\varphi):=\frac{1}{2n+1}\ln\Xi^n(\varphi)\,.
\end{equation}
The pressure  $p(\varphi)$ is defined as
\begin{equation}\label{eq10}
p(\varphi)=\lim_{n\ra\infty}P_{E^n}(\varphi)\,.
\end{equation}
The result in \eqref{eq10} is independent of the choice of the sets $E^n$.  From now on we  choose  $E^n$ so that if $x\in E^n$, then
$x_j=0$, for all $|j|>n$.

Let $[k, \ell]\subset [-n,n]$, $x\in \Sigma^\beta$. Set 
$$
x^-_k:=x_{(-\infty,k-1]}\quad \text{and}\quad x^+_\ell:=x_{[\ell+1,\infty)}\,.
$$
By our choice of $E^n$ we can
extend \eqref{eqs(v)} and \eqref{eq4} to the infinite sequence $x^-_k$ since $x_j=0$ for $j<-n$. Let
$s(x_k^-)$ be  the largest suffix $\in\cP^\beta$ of $x^-_k$, and 
$$
\widehat{x_k^-}:=y\,\widehat{s(x_k^-)}\quad\text{where $x^-_k=ys(x_k^-)$.}
$$
Let
$$
E^n_{[k,\ell]}(v):=\{x\in E^n\colon x_{[k,\ell]}=v\}\quad\text{and}\quad
E^{*,n}_{[k,\ell]}(v):=\{x\in E^n_{[k,\ell]}(v)\colon s(x_k^-)=\epsilon\}\,,
$$
and set
\begin{equation}\label{eq11}
\Xi^n_{[k,\ell]}(v)\!\!:=\!\!\sum_{x\in E^n_{[k,\ell]}(v)}\exp\big(\!\!\sum_{j\in[-n,n]}\varphi(T^jx)\big)\,,\quad
\Xi^{*,n}_{[k,\ell]}(v)\!\!:=\!\!\sum_{x\in E^{*,n}_{[k,\ell]}(v)}\!\!\exp\big(\!\!\sum_{j\in[-n,n]}\varphi(T^jx)\big)\,.
\end{equation}
We have
$$
\Xi^n(\varphi)=\sum_v\,\Xi^n_{[k,\ell]}(v)\,.
$$
Lemmas \ref{lem3} and \ref{lem4} give basic estimates used in the proof of theorem \ref{thmmain}.

\begin{lem}\label{lem3}
Let  $[k, \ell]\subset [-n,n]$ and $\|\varphi\|_\delta<\infty$. Then
$$
\Xi^{*,n}_{[k,\ell]}(v)\leq\Xi^n_{[k,\ell]}(v)\leq(z^\beta(c_1)+2){\rm e}^{2\|\varphi\|_\delta}\, \Xi^{*,n}_{[k,\ell]}(\widehat{v})\,.
$$
\end{lem}

\prf
This first inequality follows from $E^{*,n}_{[k,\ell]}(v)\subset E^n_{[k,\ell]}(v)$.
We define a map $\widehat{f}$ from $E^n_{[k,\ell]}(v)$ to $E^{*,n}_{[k,\ell]}(\widehat{v})$ by setting
\begin{equation}\label{map}
\widehat{f}:E^n_{[k,\ell]}(v)\ra E^{*,n}_{[k,\ell]}(\widehat{v})\,,\quad
x=x_k^-vx_\ell^+\mapsto  \widehat{f}(x):=\widehat{x_k^-}\widehat{v}x_\ell^+\,.
\end{equation}
Since  $\widehat{x_k^-}$and $\widehat{v}$ are presented by paths in $\cG^\beta$ with end-point $q_0$, $\widehat{f}(x)$ is well-defined and
$\widehat{f}(x)\in E^{*,n}_{[k,\ell]}(\widehat{v})$.
By lemma \ref{lem1} the map $x_k^-\mapsto \widehat{x_k^-}$ 
is at most $(z^\beta(c_1)+2)$-to-$1$.  Hence the map $\widehat{f}$ is at most $(z^\beta(c_1)+2)$-to-$1$ ($v$ is fixed).
The sequences $x$ and $x^\prime=\widehat{f}(x)$ differ  at most at two coordinates, so that by lemma \ref{lemfund}
$$
\Big|\sum_{j\in[-n,n]}\big(\varphi(T^jx)-\varphi(T^jx^\prime)\big)\Big|\leq 2 \|\varphi\|_\delta\,.
$$
Hence
\begin{eqnarray*}
\sum_{x\in E^n_{[k,\ell]}(v)}\exp\big(\sum_{j\in[-n,n]}\varphi(T^jx)\big)
&\leq&
{\rm e}^{2\|\varphi\|_\delta}
\sum_{x\in E^n_{[k,\ell]}(v)}\exp\big(\sum_{j\in[-n,n]}\varphi(T^jx^\prime)\big)\\
&\leq&
(z^\beta(c_1)+2){\rm e}^{2\|\varphi\|_\delta}
\sum_{\substack{y\in E^n_{[k,\ell]}(\widehat{v}):\\s(y^-_k)=\epsilon}}\exp\big(\sum_{j\in[-n,n]}\varphi(T^j y)\big)\\
&=&
(z^\beta(c_1)+2){\rm e}^{2\|\varphi\|_\delta} \,\Xi^{*,n}_{[k,\ell]}(\widehat{v})\,.
\end{eqnarray*}
\qed

\begin{lem}\label{lem4}
Let  $[k, \ell]\subset [-n,n]$ and $\|\varphi\|_\delta<\infty$. 
Then 
$$
\Xi^{*,n}_{[k,\ell]}(v)\geq |\tA|^{-(z^\beta(v)+1)}{\rm e}^{-(z^\beta(v)+2)\|\varphi\|_\delta}\,\Xi^{*,n}_{[k,\ell]}(\widehat{v})\,.
$$
\end{lem}

\prf
If $s(v)=\epsilon$, then $v=\widehat{v}$ and the inequality is trivial.
Let $v$ with $s(v)\not=\epsilon$; we define a map $f\colon E^{*,n}_{[k,\ell]}(\widehat{v})\ra E^{*,n}_{[k,\ell]}(v)$,
$$
f(x_k^-\widehat{v}x_\ell^+):=x_k^- v{x^\prime}_\ell^+\,,
$$
with
$$
{x^\prime}_{\ell+1}^+=\cdots={x^\prime}_{\ell+z^\beta(v)+1}^+=0\quad\text{and}\quad {x^\prime}_{j}^+=x_j\quad\text{if $j>\ell+z^\beta(v)+1$.}
$$
This map is at most $|\tA|^{z^\beta(v)+1}$-to-$1$. 
\begin{eqnarray*}
\sum_{j\in [-n,n]}\varphi(T^j(x_k^-\widehat{v}x_\ell^+))
&=&
\sum_{j\in [-n,n]}\big(\varphi(T^j(x_k^-\widehat{v}x_\ell^+))-\varphi(T^j(x_k^-\widehat{v}{x^\prime}_\ell^+))\big)\\
&+&
\sum_{j\in [-n,n]}\big(\varphi(T^j(x_k^-\widehat{v}{x^\prime}_\ell^+))-\sum_{j\in [-n,n]}\varphi(T^j(x_k^-v{x^\prime}_\ell^+))\big)\\
&+&
\sum_{j\in [-n,n]}\varphi(T^j(x_k^-v{x^\prime}_\ell^+))\,.
\end{eqnarray*}
The configurations $x_k^-\widehat{v}x_\ell^+$ and $x_k^-\widehat{v}{x^\prime}_\ell^+$ differ  at most at 
$i=\ell+1,\ldots,\ell+z^\beta(v)+1$, so that
\begin{equation}\label{delta}
\Big|\sum_{j\in [-n,n]}\big(\varphi(T^j(x_k^-\widehat{v}x_\ell^+))-\varphi(T^j(x_k^-\widehat{v}{x^\prime}_\ell^+))\big)\Big|
\leq (z^\beta(v)+1)\|\varphi\|_\delta\,.
\end{equation}
The configurations $x_k^-\widehat{v}{x^\prime}_\ell^+$ and $x_k^-v{x^\prime}_\ell^+$ differs at one coordinate, so that
$$
\Big|\sum_{j\in [-n,n]}\big(\varphi(T^j(x_k^-\widehat{v}{x^\prime}_\ell^+))-
\sum_{j\in [-n,n]}\varphi(T^j(x_k^-v{x^\prime}_\ell^+))\big)\Big|\leq\|\varphi\|_\delta\,.
$$
Therefore
\begin{eqnarray*}
\Xi^{*,n}_{[k,\ell]}(\widehat{v})
&=&
\sum_{x\in E^{*,n}_{[k,\ell]}(\widehat{v})}\exp\sum_{j\in [-n,n]}\varphi(T^jx)\\
&=&
\sum_{y\in E^{*,n}_{[k,\ell]}(v)}\sum_{\substack{x\in E^{*,n}_{[k,\ell]}(\widehat{v}):\\f(x)=y}}
\exp\sum_{j\in [-n,n]}\varphi(T^jx)\\
&\leq&
|\tA|^{z^\beta(v)+1}{\rm e}^{(z^\beta(v)+2)\|\varphi\|_\delta}
\sum_{y\in E^{*,n}_{[k,\ell]}(v)}\exp\sum_{j\in [-n,n]}\varphi(T^jy)\\
&=&
|\tA|^{z^\beta(v)+1}{\rm e}^{(z^\beta(v)+2)\|\varphi\|_\delta}\,\Xi^{*,n}_{[k,\ell]}(v)\,.
\end{eqnarray*}
\qed

\begin{lem}\label{lem5}
Let $w\in\cL^\beta_m$ and $w^\sharp\in\Sigma^\beta$,
\begin{equation}\label{sharp}
w^\sharp_j:=\begin{cases}
w_j&\text{if $1\leq j\leq m$}\\
0 &\text{otherwise.}\\
\end{cases}
\end{equation}
The pressure $p(\varphi)$ is equal to
$$
\lim_{m\ra\infty}\frac{1}{m}\ln\!\!\!\!\sum_{w\in\cL^\beta_m:s(w)=\epsilon}\!\!\!\!\exp\sum_{j=1}^m\varphi(T^jw^\sharp)\,.
$$
\end{lem}
\prf
The configurations $w^\sharp$ and ${\widehat{w}}^\sharp$ differ  at most at one coordinate, say coordinate $i$.
Therefore
$$
\Big|\sum_{j=1}^m\big(\varphi(T^j{\widehat{w}}^\sharp)-\varphi(T^jw^\sharp)\big)\Big|
\leq\|\varphi\|_\delta\,.
$$
The map $w\mapsto\widehat{w}$ is at most $(z^\beta(c_1)+2)$-to-$1$ (lemma \ref{lem1}).
\qed

\subsection{Main result}\label{subsectionresult}

\begin{defn}\label{defnweakgibbs}
An invariant probability measure $\nu$ is  a \emph{weak Gibbs measure for a continuous function $\psi$},
if for any $\delta>0$ there exists $N_\delta$ such that for $m\geq N_\delta$,
$$
\sup_{x\in X}\Big|\frac{1}{m}\ln\nu([x_0\cdots x_{m-1}])-\frac{1}{m}\sum_{\ell=0}^{m-1}\psi(T^\ell(x))\Big|\leq \delta\,,
$$
where $[x_0\cdots x_{m-1}]=\{y\in\Sigma^\beta\colon y_0\cdots y_{m-1}=x_0\cdots x_{m-1}\}$.
\end{defn}

\noindent
{\bf Remark.\,}  Our definition of weak Gibbs measure is stated so that if $\nu$ is weak Gibbs for $\psi$, then $p(\psi)=0$ \cite{PS2}. If 
$\psi=\varphi-p(\varphi)$, then it equivalent to
$$
{\rm e}^{-\delta m}\leq
\frac{\nu([x_0\cdots x_{m-1}])}{\exp\big(-mp(\varphi)+\sum_{\ell=0}^{m-1}\varphi(T^\ell x)\big)}\leq {\rm e}^{\delta m}\,.
$$
\qed

\begin{thm}\label{thmmain}
Let $\beta>1$ and  $\varphi$ be a function of bounded total oscillations on $\Sigma^\beta$. 

\noindent
1) If $\nu$ is an equilibrium measure for $\varphi$ and if
$$
\lim_{n\ra\infty}\frac{\overline{z}^\beta(n)}{n}=0\,,
$$
then $\nu$ is a weak Gibbs measure for $\psi=\varphi-p(\varphi)$.

\noindent
2)
If $\nu$ is an equilibrium measure for $\varphi$ and if
$$
\limsup_{n\ra\infty}\frac{\overline{z}^\beta(n)}{n}>0\,,
$$
then $\nu$ is not a weak Gibbs measure for $\psi=\varphi-p(\varphi)$.
\end{thm}
If $\nu$ is a weak Gibbs measure, then the empirical measures verify a large deviations principle \cite{PS2}. 
Large deviations for (one-sided) $\beta$-shifts, for any $\beta>1$ and equilibrium measures have been proved 
by Climenhaga, Thompson and Yamamoto \cite{CTY} for the class of functions verifying Bowen condition. From the estimates of  lemma \ref{lem7}, proposition 4.3 and theorem 3.1 in \cite{PS1} the result of \cite{CTY} is also valid for all equilibrium measures for functions $\varphi$ of bounded total oscillations. This is important since Bowen condition implies uniqueness of the equilibrium measure for $\beta$-shifts, while this is not necessarily the case for bounded total oscillations functions. 

\section{Proof of theorem \ref{thmmain}}\label{sectionproof}

Let $\varphi$ be a function of bounded total oscillations on $\Sigma^\beta$. In subsection \ref{subsection3.1} we prove upper and lower bound for 
$\nu([y_0\cdots y_{m-1}])$
for any equilibrium measure $\nu$ of $\varphi$. There is no restriction on $\beta>1$. In subsection \ref{subsection3.2}
we prove theorem \ref{thmmain}.

\subsection{Upper and lower bounds}\label{subsection3.1}

We first assume  that there is a unique tangent functional $\nu$ to the pressure at $\varphi$.
The result  is then extended to any $\varphi$ of bounded total oscillations 
using a theorem of Mazur and  a theorem of Lanford and Robinson (see e.g. \cite{Ru} appendix A.3.7).

When there is a unique tangent functional to the pressure at $\varphi$ we can estimate
$\nu([y_0\cdots y_{m-1}])$ using a classical result about
differentiability of a convex function, here the pressure, which is a pointwise limit of convex functions, theorem 25.7 in \cite{Ro}.
Let $\bar u\in\cL^\beta_m$ be fixed and set
\begin{equation*}
I_{\bar u}(y):=\left\{\begin{array}{lll}
1 &\text{if $y_0\cdots y_{m-1}=\bar u$}\\
0 &\text{otherwise.}
\end{array}\right.
\end{equation*}
We have
\begin{eqnarray}\label{eq13}
  \nu(I_{\bar u})
  &=&
  \left .\frac {d}{dt} \lim_{n\to\infty}
  P_{E^n}(\varphi+t\,I_{\bar u})\right|_{t=0}=
 \lim_{n\to\infty}\left .\frac {d}{dt} P_{E^n}(\varphi+t\,I_{\bar u})\right|_{t=0}\nonumber\\
 &=&
\lim_{n\ra\infty} \frac{1}{2n+1}\sum_{j=-n}^n
 \frac
{ \sum_{x\in E^n}I_{\bar u}(T^j x) \exp\sum_{i=-n}^n \varphi(T^i x)}
{ \sum_{x\in E^n}\exp\sum_{i=-n}^n \varphi(T^i x)}\,.
\end{eqnarray}
Let $j=k$ and $\ell:=k+m-1$. Then a term in \eqref{eq13} is written as a ratio of partition functions (see \eqref{eq11})
\begin{equation}\label{eq14}
 \frac
{ \sum_{x\in E^n}I_{\bar u}(T^j x) \exp\sum_{i=-n}^n \varphi(T^i x)}
{ \sum_{x\in E^n}\exp\sum_{i=-n}^n \varphi(T^i x)}=
\frac{\Xi^{n}_{[k,\ell]}(\bar u)}{\Xi^n(\varphi)}=\frac{\Xi^{n}_{[k,\ell]}(\bar u)}{\sum_v \Xi^{n}_{[k,\ell]}(v)}\,.
\end{equation}
The core of the proof involves estimating the ratio of partition functions $\Xi^{n}_{[k,\ell]}(v)/\Xi^{n}_{[k,\ell]}(\bar u)$, uniformly in 
$[k,\ell]\subset[-n,n]$
using lemmas \ref{lem3} and \ref{lem4}. Since in \eqref{eq13} we take the limit $n\ra\infty$, it is sufficient to consider
the cases where $[k,\ell]\subset[-n,n]$.

\begin{lem}\label{lem6}
Let $\bar u\in\cL^\beta_m$. For any $\varepsilon>0$ and continuous $\varphi$ such that $\|\varphi\|_\delta<\infty$ there exists 
$N_{\varepsilon,\varphi}$ so that if
$y$ is such that $J_{[1,m]}(y)=\bar u$ and
$m\geq N_{\varepsilon,\varphi}$, then
$$
\nu(I_{\bar u})\leq  K_{\varphi,\varepsilon}^+(m,\beta)  \exp\big(\sum_{j=1}^m\varphi(T^jy)-mp(\varphi)\big)\,,
$$
where 
$$
 K_{\varphi,\varepsilon}^+(m,\beta)=(z^\beta(c_1)+2){\rm e}^{3\|\varphi\|_\delta}{\rm e}^{5m\varepsilon}\,.
$$
$N_{\varepsilon,\varphi}$ is chosen  so that all of the inequalities \eqref{eq6}, \eqref{eq7} and \eqref{N} are satisfied for 
$m\geq N_{\varepsilon,\varphi}$.
\end{lem}

\prf
Let $\varepsilon>0$ and $\varphi$ be given and 
$N_{\varepsilon,\varphi}$  defined as above. 
We consider the term in  the sum \eqref{eq13} with  $j=k$ and
$[k,\ell]\subset[-n,n]$ (see \eqref{eq14}),
$$
\frac{\Xi^{n}_{[k,\ell]}(\bar u)}{\Xi^n(\varphi)}=\frac{\Xi^{n}_{[k,\ell]}(\bar u)}{\sum_v \Xi^{n}_{[k,\ell]}(v)}\,.
$$
We have
$$
\sum_v \Xi^{n}_{[k,\ell]}(v)\geq \sum_{v:s(v)=\epsilon} \Xi^{n}_{[k,\ell]}(v)\geq  \sum_{v:s(v)=\epsilon} \Xi^{*,n}_{[k,\ell]}(v)\,.
$$
From now on $s(v)=\epsilon$. By lemma \ref{lem3}
$$
\Xi^n_{[k,\ell]}(\bar u)\leq(z^\beta(c_1)+2){\rm e}^{2\|\varphi\|_\delta}\, \Xi^{*,n}_{[k,\ell]}(\widehat{\bar u})\,.
$$
Hence
$$
\frac{\Xi^{*,n}_{[k,\ell]}(v)}{\Xi^n_{[k,\ell]}(\bar u)}\geq \frac{1}{(z^\beta(c_1)+2){\rm e}^{2\|\varphi\|_\delta}}\,
 \frac{\Xi^{*,n}_{[k,\ell]}(v)}{\Xi^{*,n}_{[k,\ell]}(\widehat{\bar u})}\,.
$$
Since $s(v)=\epsilon$ and $s(\widehat{\bar u})=\epsilon$, if $s(x_k^-)=\epsilon$, then 
\begin{equation}\label{eq15}
x_k^-vx_\ell^+\in E^{*,n}_{[k,\ell]}(v) \iff x_k^-\widehat{\bar u}x_\ell^+\in E^{*,n}_{[k,\ell]}(\widehat{\bar u})\,.
\end{equation}
We write
\begin{equation}\label{eq16}
\Xi^{*,n}_{[k,\ell]}(v) 
=\sum_{x\in E^{*,n}_{[k,\ell]}(v)}
\frac{\exp\big(\sum_{j=-n}^n\varphi(T^j( x_k^-vx_\ell^+))}{\exp\big(\sum_{j=-n}^n\varphi(T^j( x_k^-\widehat{\bar u}x_\ell^+))\big)}
\exp\big(\sum_{j=-n}^n\varphi(T^j( x_k^-\widehat{\bar u}x_\ell^+))\big)\,.
\end{equation}
Let  $v^\sharp$ and ${\widehat{\bar u}}^\sharp$ be defined as in \eqref{sharp} with $1\leq j\leq m$ replaced by $k\leq j\leq \ell$.
\begin{eqnarray}\label{ineq1}
& &\sum_{j=-n}^n\big(\varphi(T^j( x_k^-vx_\ell^+))-\varphi(T^j( x_k^-\widehat{\bar u}x_\ell^+))\big)
=
\sum_{j\not\in[k,\ell]}\big(\varphi(T^j( x_k^-vx_\ell^+))-\varphi(T^j( x_k^-\widehat{\bar u}x_\ell^+))\big) +\nonumber\\
& &\sum_{j\in[k,\ell]}\big(\varphi(T^j( x_k^-vx_\ell^+))-\varphi(T^jv^\sharp)\big)+
\sum_{j\in[k,\ell]}\big(\varphi(T^j{\widehat{\bar u}}^\sharp)-\varphi(T^j( x_k^-\widehat{\bar u}x_\ell^+))\big)
+\\
& &\sum_{j\in[k,\ell]}\big(\varphi(T^jv^\sharp)-\varphi(T^j{\widehat{\bar u}}^\sharp)\big)\,.\nonumber
\end{eqnarray}
By lemma \ref{lem2}, if $m\geq N_{\varepsilon,\varphi}$, then
$$
\Big|\sum_{j=-n}^n\big(\varphi(T^j( x_k^-vx_\ell^+))-\varphi(T^j( x_k^-\widehat{\bar u}x_\ell^+))\big)
-\sum_{j\in[k,\ell]}\big(\varphi(T^jv^\sharp)-\varphi(T^j{\widehat{\bar u}}^\sharp)\big)\Big|\leq 3\varepsilon m\,.
$$
Let $y\in\Sigma^\beta$ be such that $J_{[k,\ell]}(y)=\bar u$. By definition  of  ${\bar u}^\sharp$ 
$$
\text{$J_{[k,\ell]}({\bar u}^\sharp)=\bar u$ and  ${\bar u}^\sharp_i=0$ for all $i\not\in [k,\ell]$.}
$$
By lemmas \ref{lemfund} and \ref{lem2}, if $m\geq N_{\varepsilon,\varphi}$, then
\begin{eqnarray}\label{ineq2}
& &\Big|\sum_{j\in[k,\ell]}\big(\varphi(T^j{\widehat{\bar u}}^\sharp)-\varphi(T^jy)\big)\Big|
\leq \\ 
& &\Big|\sum_{j\in[k,\ell]}\big(\varphi(T^j{\widehat{\bar u}}^\sharp)-\varphi(T^j{\bar u}^\sharp)\big)\Big|+
\Big|\sum_{j\in[k,\ell]}\big(\varphi(T^j{\bar u}^\sharp)-\varphi(T^jy)\big)\Big|\leq\|\varphi\|_\delta +m\varepsilon\,.\nonumber
\end{eqnarray}
$N_{\varepsilon,\varphi}$ has been chosen so that (see lemma \ref{lem5})
\begin{equation}\label{N}
{\rm e}^{-m\varepsilon}\leq\frac{{\rm e}^{mp(\varphi)}}{\sum_{\substack{v\in\cL^\beta_m:\\
s(v)=\epsilon}}\exp\sum_{j=1}^m\varphi(T^jv^\sharp)}
\leq{\rm e}^{m\varepsilon}\,.
\end{equation}
From   \eqref{eq14}, \eqref{eq16}, $s(v)=\epsilon$ and the above estimates, taking into account \eqref{eq15}, which allows the use of the 
elementary inequalities for positive real numbers
$a_i$ and $b_i$,
\begin{equation}\label{elementary}
\inf_i \frac{a_i}{b_i}\leq\frac{\sum_{i=1}^n a_i}{\sum_{i=1}^nb_i }=\frac{\sum_{i=1}^n \frac{a_i}{ b_i} \,b_i} {\sum_{i=1}^n b_i}\leq
\sup_i \frac{a_i}{b_i}\,,
\end{equation}
we get
\begin{eqnarray}\label{eq19}
\frac{\sum_v \Xi^{n}_{[k,\ell]}(v)}{\Xi^{n}_{[k,\ell]}(\bar u)}
&\geq &
\frac{\sum_{v:s(v)=\epsilon} \Xi^{*,n}_{[k,\ell]}(v)}{\Xi^{n}_{[k,\ell]}(\bar u)}
\geq
\frac{1}{(z^\beta(c_1)+2){\rm e}^{2\|\varphi\|_\delta}}\,
\frac{\sum_{v:s(v)=\epsilon} \Xi^{*,n}_{[k,\ell]}(v)}{\Xi^{n}_{[k,\ell]}(\widehat{\bar u})}\nonumber\\
&\geq &
{\rm e}^{-4m\varepsilon-\|\varphi\|_\delta}\frac{{\rm e}^{-\sum_{j\in[k,\ell]}\varphi(T^jy)}}{(z^\beta(c_1)+2){\rm e}^{2\|\varphi\|_\delta}}
\sum_{v:s(v)=\epsilon}\exp\sum_{j=1}^m\varphi(T^jv^\sharp)\nonumber\\
&\geq&
\frac{{\rm e}^{-5m\varepsilon}}{(z^\beta(c_1)+2){\rm e}^{3\|\varphi\|_\delta}}
\exp\big(-\sum_{j\in[k,\ell]}\varphi(T^jy)+mp(\varphi)\big)\,.\nonumber
\end{eqnarray}
The result follows from \eqref{eq13} by taking the limit $n\ra\infty$.
\qed

\begin{lem}\label{lem7}
Let $\bar u\in\cL^\beta_m$. For any $\varepsilon>0$ and continuous $\varphi$ such that $\|\varphi\|_\delta<\infty$ there exists 
$N_{\varepsilon,\varphi}$ so that if
$y$ is such that $J_{[1,m]}(y)=\bar u$ and
$m\geq N_{\varepsilon,\varphi}$, then
$$
\nu(I_{\bar u})\geq
K^-_{\varphi,\varepsilon}(m,\beta,\bar{u})
\exp\big(\sum_{j=1}^m\varphi(T^jy)-mp(\varphi)\big)\,,
$$
where 
$$
K^-_{\varphi,\varepsilon}(m,\beta,\bar{u})=
\frac{ |\tA|^{-(z^\beta(\bar u)+1)}{\rm e}^{-(z^\beta(\bar u)+2)\|\varphi\|_\delta}}{(z^\beta(c_1)+2)^2 {\rm e}^{5m\varepsilon +3\|\varphi\|_\delta}}\,.
$$
$N_{\varepsilon,\varphi}$ is chosen  so that all of the inequalities \eqref{eq6}, \eqref{eq7} and \eqref{N} are satisfied for 
$m\geq N_{\varepsilon,\varphi}$.
\end{lem}

\prf
Let $\varepsilon>0$ and $\varphi$ be given and 
$N_{\varepsilon,\varphi}$  defined as above.
We consider the term in  the sum \eqref{eq13} with  $j=k$ and
$[k,\ell]\subset[-n,n]$ (see \eqref{eq14}),
$$
\frac{\Xi^{n}_{[k,\ell]}(\bar u)}{\Xi^n(\varphi)}=\frac{\Xi^{n}_{[k,\ell]}(\bar u)}{\sum_v \Xi^{n}_{[k,\ell]}(v)}\,.
$$
As in the proof of lemma \ref{lem6}, if
$s(v)=\epsilon$, then  we can estimate the ratio
$\Xi^{*,n}_{[k,\ell]}(v)/\Xi^{*,n}_{[k,\ell]}(\widehat{\bar u})$ using  \eqref{elementary}
since \eqref{eq15} holds. 
If $y\in\Sigma^\beta$ be such that $J_{[k,\ell]}(y)=\bar u$ and $s(v)=\epsilon$, then (see \eqref{ineq1} and \eqref{ineq2})
$$
\frac{\Xi^{*,n}_{[k,\ell]}(v)}{\Xi^{*,n}_{[k,\ell]}(\widehat{\bar u})}\leq 
{\rm e}^{4m\varepsilon +\|\varphi\|_\delta}\exp\big(\sum_{j\in[k,\ell]}\varphi(T^jv^\sharp)-\varphi(T^jy)\big)\,.
$$
Let $m\geq N_{\varepsilon,\varphi}$.
\begin{eqnarray*}
\frac{\sum_v \Xi^{n}_{[k,\ell]}(v)}{\Xi^{n}_{[k,\ell]}(\bar u)}
&\leq & 
\sum_v \frac{(z^\beta(c_1)+2){\rm e}^{2\|\varphi\|_\delta}}{ |\tA|^{-(z^\beta(\bar u)+1)}{\rm e}^{-(z^\beta(\bar u)+2)\|\varphi\|_\delta}}\,
\frac{\Xi^{*,n}_{[k,\ell]}(\widehat{v})}{\Xi^{*,n}_{[k,\ell]}(\widehat{\bar u})} \\
&\leq&
 \frac{(z^\beta(c_1)+2)^2{\rm e}^{2\|\varphi\|_\delta}}{ |\tA|^{-(z^\beta(\bar u)+1)}{\rm e}^{-(z^\beta(\bar u)+2)\|\varphi\|_\delta}}\,
\sum_{v:s(v)=\epsilon}\frac{\Xi^{*,n}_{[k,\ell]}(v)}{\Xi^{*,n}_{[k,\ell]}(\widehat{\bar u})} \\
&\leq &
 \frac{(z^\beta(c_1)+2)^2 {\rm e}^{5m\varepsilon +3\|\varphi\|_\delta}}{ |\tA|^{-(z^\beta(\bar u)+1)}{\rm e}^{-(z^\beta(\bar u)+2)\|\varphi\|_\delta}}\,
 \exp\big(\sum_{j\in[k,\ell]}-\varphi(T^jy)+mp(\varphi)\big)\,.
\end{eqnarray*}
For the first inequality we use lemmas \ref{lem3}, \ref{lem4}, and for the second inequality, where we replace $\widehat{v}$ by $v$ with $s(v)=\epsilon$,
we use lemma \ref{lem1} and $s(\widehat{v})=s(v)=\epsilon$.
The result follows from \eqref{eq13} by taking the limit $n\ra\infty$.
\qed

We now remove the restriction that $\nu$ is the unique equilibrium measure for $\varphi$.
The pressure is convex and continuous on the Banach space $\cB(\Sigma^\beta)$ of bounded total oscillations functions.
The set $\cR\subset\cB$ of $\varphi$ such that $\partial p(\varphi)=\{\nu\}$ has a unique element $\nu$ is residual
(theorem of Mazur). 

Let $\varphi\in\cB$ be an arbitrary function of bounded total oscillations and $\nu$ be an equilibrium measure for $\varphi$, such that there exists
a sequence  $\varphi_k\in\cR$ with the properties that $\lim_k\varphi_k=\varphi$ and $\lim_k\nu_k=\nu$, $\{\nu_k\}=\partial p(\varphi_k)$
 (weak convergence). For that $\varphi$ let $N_{\varepsilon,\varphi}$ be defined as in lemmas \ref{lem6} and \ref{lem7}. Let
 $\varepsilon^\prime>\varepsilon$. By our choice of $r_\varepsilon$ (see proof of lemma \ref{lem2})
 $$
 \sum_{k:|k|>r_\varepsilon}\delta_k(\varphi_k)\leq \sum_{k:|k|>r_\varepsilon}\delta_k(\varphi)+\|\varphi_k-\varphi\|_\delta
 \leq \frac{\varepsilon}{2}+ \|\varphi_k-\varphi\|_\delta\,,
 $$
and $\|\varphi_k\|\leq\|\varphi\|+\|\varphi_k-\varphi\|$. Since $|p(\varphi)-p(\varphi_k)|\leq\|\varphi-\varphi_k\|$,
for $m\geq N_{\varepsilon,\varphi}$,
 the upper and lower bounds of lemmas \ref{lem6} and \ref{lem7} are true for $\nu_k$ and $\varphi_k$, with constants
 $K_{\varphi_k,\varepsilon^\prime}^+(m,\beta)$ and  $K^-_{\varphi_k,\varepsilon^\prime}(m,\beta,\bar{u})$,
provided that $k$ is large enough. Since $\varepsilon^\prime>\varepsilon$ is arbitrary,
 lemmas \ref{lem6} and \ref{lem7} are true  for $\nu$.
 This is also the case for any $\mu$ in the weak-closed convex hull of such $\nu$'s. By the theorem of Lanford and Robinson this set coincides with the set of equilibrium measures for $\varphi$.
\qed

\subsection{Proof of theorem \ref{thmmain}}\label{subsection3.2}

1) Suppose that $\lim_{n\ra\infty}\overline{z}^\beta(n)/n=0$. Then 
$$
\lim_{m\ra 0}\sup_{\bar u\in\cL^\beta_m}\frac{z^\beta(\bar u)}{m}=0\,.
$$
Tthe estimates of lemmas \ref{lem6} and \ref{lem7} prove that the equilibrium measure $\nu$ for $\varphi$ is
a weak Gibbs measure for $\psi=\varphi-p(\varphi)$.

 \noindent
2) Suppose that $\limsup_{n\ra\infty}\overline{z}^\beta(n)/n>0$. There exists
an increasing diverging sequence $\{m_k\}_k$ and $w^k\in\cL^\beta_{m_k}$ so that
$\lim_kz^\beta(w^k)/m_k=a>0$. 
Let
$$
\widetilde {w^k}:=w^k\underbrace{0\cdots 0}_{z^\beta(w^k)}\,.
$$
By definition $|\widetilde{w_k}|=m_k+z^\beta(w^k)$ and $\nu(I_{w^k})=\nu(I_{\widetilde{w^k}})$.
If $y\in[w^k]$, then by lemma \ref{lem6} (and $m_k$ large enough)
\begin{eqnarray*}
\frac{\nu(I_{w^k})}{\exp\big(-m_kp(\varphi)+\sum_{\ell=1}^{m_k}\varphi(T^\ell y)\big)}
=
\frac{\nu(I_{\widetilde{w^k}})}{\exp\big(-m_kp(\varphi)+\sum_{\ell=1}^{m_k}\varphi(T^\ell y)\big)}\leq\\
K^+_{\varphi,\varepsilon}(m_k+z^\beta(w^k),\beta)\exp\Big(\sum_{j=1}^{z^\beta(w^k)}(\varphi(T^{m_k+j}y)-p(\varphi))\Big)\,.
\end{eqnarray*}
We can compare $\sum_{j=1}^{z^\beta(w^k)}\varphi(T^{m_k+j}y)$ with 
$z^\beta(w^k)\varphi(0)$, where $0$ is the configuration with all coordinates equal to $0$. For $z^\beta(w^k)\geq N_{\varepsilon,\varphi}$,
$$
\Big|\sum_{j=1}^{z^\beta(w^k)}(\varphi(T^{m_k+j}y))-z^\beta(w^k)\varphi(0)
\Big|\leq \varepsilon z^\beta(w^k)\,,
$$
so that 
$$
\frac{\nu(I_{w^k})}{\exp\big(-m_kp(\varphi)+\sum_{\ell=1}^{m_k}\varphi(T^\ell y )\big)}\leq
K^+_{\varphi,\varepsilon}(m_k+z^\beta(w^k),\beta){\rm e}^{z^\beta(w^k)(\varphi(0)-p(\varphi)+\varepsilon)}\,.
$$
For any $T$-invariant probability measure $\mu$ and any continuous function $\varphi$,
$$
p(\varphi)\geq h_T(\mu)+\int\varphi\,d\mu\,,
$$ 
where $h_T(\mu)$ is the (metric) entropy of $\mu$. A $T$-invariant probability measure $\nu$ is an equilibrium measure for $\varphi$ if and only if 
$$
p(\varphi)=h_T(\nu)+\int\varphi\,d\nu\,.
$$
The support of a measure is the complement of the union of the open sets of measure zero.
For the $\beta$-shift any equilibrium measure has support $\Sigma^\beta$, since by lemma \ref{lem7} all cylinder sets have positive measure and the cylinder sets generate the topology. Therefore the Dirac measure $\delta_0$ cannot be an equilibrium measure. Hence
\begin{equation}\label{eqp}
p(\varphi)>h_T(\delta_0)+\int \varphi\,\delta_0=\varphi(0)
\end{equation}
because the entropy of $\delta_0$ is zero.
Inequality \eqref{eqp} implies that $\nu$ is not a weak Gibbs measure for $\psi=\varphi-p(\varphi)$.
Indeed,
for $x$ such that $J_{[1,m_k]}(x)=w_k$, 
\begin{eqnarray*}
\frac{1}{m_k}\Big(\ln\nu (I_{w^k})-\sum_{j=1}^{m_k}\psi(T^jx)\Big)
&\leq&
\frac{\ln K^+_{\varphi,\varepsilon}(m_k+z^\beta(w^k),\beta)}{m_k}\\
&+&
\frac{z^\beta(w^k)(\varphi(0)-p(\varphi)+\varepsilon)}{m_k}\,.
\end{eqnarray*}
Taking the limit $k\ra\infty$ and observing that   $\varepsilon$ is as small as we wish in that limit,
$$
\lim_{k\ra\infty}\frac{1}{m_k}\Big(\ln\nu (I_{w^k})-\sum_{j=1}^{m_k}\psi(T^jx)\Big)
\leq
\lim_{k\ra\infty}\frac{z^\beta(w^k)}{m_k}\underbrace{\big(6\varepsilon+\varphi(0)-p(\varphi)\big)}_{\text{$<0$ if $\varepsilon$ is small enough}}
<0\,.
$$
\qed

\end{document}